\documentstyle[11pt,leqno,graphicx,amscd,fancyhdr,amssymb,latexsym,amsbsy,hyperref,xypic,mathrsfs,pstricks,color,verbatim]{amsart}

\def\da{\delta}
\def\gz{\gamma}
\def\ggz{\Gamma}

\def\ra{\rightarrow}

\def\t{\tau}

\def\Int{\mbox{Int}}

\newtheorem{thm}{Theorem}[section]
\newtheorem{lem}[thm]{Lemma}

\newtheorem{prop}[thm]{Proposition}
\newtheorem{rem}[thm]{Remark}

\theoremstyle{definition}
\newtheorem{example}[thm]{Example}

\newtheorem{defn}[thm]{Definition}

\begin{document}

\title{Non-leaving-face property for marked surfaces}
\author{Thomas Br\"ustle and Jie Zhang}


\begin{abstract}We consider the polytope arising from a marked surface by flips of triangulations. Sleator, Tarjan and Thurston studied in 1988 the diameter of the associahedron, which is the polytope arising from a marked disc by flips of triangulations. They showed that every shortest path between two vertices in a face does not leave that face.
We establish that same non-leaving-face property for all unpunctured marked surfaces.
\end{abstract}
\maketitle


\section{Introduction}

The exchange graph is a central notion in the theory of cluster algebras initiated by S. Fomin
and A. Zelevinsky in \cite{Fomin-Zelevinsky2002-2,Fomin-Zelevinsky2003-2} and in the theory of cluster categories  introduced in \cite{BMRRT}. The vertices of the exchange graph are given in this context by clusters, and edges between two vertices are given by mutations of the corresponding clusters.
The exchange graph has the structure of a generalized (or abstract) polytope, see \cite{Fomin-Zelevinsky2003-2, Chapoton-Fomin-Zelevinsky2002,Reading2006,Hohlweg-Lange-Thomas2011}.
\medskip

The non-leaving-face property of a polytope was introduced in \cite{Ceballos-Pilaud2016} by C. Ceballos and V. Pilaud, and further studied in \cite{Williams2015} by N. Williams. We say a polytope $P$ has the non-leaving-face property if any shortest path connecting two vertices in the graph of $P$ stays in the minimal face of P containing both.
This property has first been established for the n-dimensional associahedron of type A in \cite{Sleator-Tarjan-Thurston1988}, with the aim to find the diameter of these associahedra.  C. Ceballos and V. Pilaud proved in \cite{Ceballos-Pilaud2016} that associahedra of types B,C and D also have the non-leaving-face property. Moreover, N. Williams  established in  \cite{Williams2015} the non-leaving-face property of W-permutahedra
and W-associahedra, for a finite Coxeter system (W,S). We like to mention that not all (generalized) associahedra  satisfy the non-leaving-face property, see \cite{Ceballos-Pilaud2016} for more details.
\smallskip

We study in this paper the non-leaving-face property of an exchange graph coming from an unpunctured  marked surface $(S,M)$, where $S$ denotes the surface and $M$ the set of marked points on the boundary of $S$.
This exchange graph, as cluster exchange graph of the cluster algebra associated with $(S,M)$, or the cluster category $(S,M)$, has been introduced in \cite{Fomin-Shapiro-Thurstion2008}, and since then intensely studied in various papers, see \cite{Labardini2009,Brustle-Zhang2011,Brustle-zhang2013,Brustle-Yu2015} and others.
However, the question of shortest paths of mutations has not been addressed in this context previously.

If the surface $S$ is a disc, then the exchange graph is an acssociahedron of type A. All other cases of unpunctured marked surfaces result in an infinite exchange graph, thus results from \cite{Sleator-Tarjan-Thurston1988}, \cite{Ceballos-Pilaud2016}, \cite{Williams2015} do not apply.
We study these infinite polytopes, arising as exchange graphs, and show that they have the non-leaving-face property (see Theorem \ref{main-theorem}).

\section{Preliminaries}

\subsection{Exchange graph and non-leaving-face property}
Following \cite{Fomin-Zelevinsky2003-2}, the definition of an exchange graph is formalized in \cite{Brustle-Yang2013} as follows:
Consider a set $V$ with a compatibility relation $R$,
that is, a reflexive and symmetric relation $R$  on $V$. We say that two elements $x$ and $y$ of $V$ are
compatible if $(x, y)\in R$.
Motivated by cluster theory, maximal subsets of pairwise compatible elements are called clusters.
Assume the following conditions:
\begin{itemize}
 \item[$(1)$] All clusters are finite
and have the same cardinality, say $n$;
\item[$(2)$] Any subset of $n-1$ pairwise compatible elements is contained in precisely
two clusters.
\end{itemize}
We then define an exchange graph to be the graph whose vertices are the clusters
and where two clusters are joined by an edge precisely when their intersection has
cardinality $n -1$. We refer to the edges of an exchange graph as mutations. Note that all exchange graphs are $n-$regular.
The conditions on the compatibility relation $R$ can be rephrased as follows:
consider the (abstract) simplicial complex $\Delta$ whose $l-$simplices are the subsets of
$l + 1$ pairwise compatible elements of $V$. A simplex of codimension 1 is called a
wall. We assume that
\begin{itemize}
  \item[$(1)$] $\Delta$ is a pure simplicial complex, i.e. all maximal simplices are of the same
dimension;
  \item[$(2)$]every wall is contained in precisely two maximal simplices.
\end{itemize}
Then the exchange graph is the dual graph of $\Delta$.
If in addition the exchange graph
is connected, then $\Delta$ is a pseudo-manifold, see \cite[Section 2.1]{Fomin-Zelevinsky2003-2}.
In this light, the faces of an exchange graph ${\mathbf{G}}$ are subgraphs corresponding to some $l-$simplex. More precisely, a face $F_U$ is the full subgraph of ${\mathbf{G}}$ given by all vertices (clusters) containing the set $U$, with $U$ being a set of pairwise compatible elements of $V$.

Exchange graphs appear originally in cluster theory, but since then many other structures have found to yield exchange graphs, such as  support-$\tau-$tilting modules over a finite-dimensional algebra, silting objects in a derived category etc,  we refer \cite{Brustle-Yang2013}
for more details.

If  two clusters $v_1,v_2$ in an exchange graph ${\mathbf{G}}$ are are joined by an edge, we denote this by
                          $$v_1-v_2$$
We call a path in ${\mathbf{G}}$ between $v$ and $w$
$$v=v_1-v_2-\ldots-v_n=w$$
a \emph{geodesic} connecting vertices $v$ to $w$ if the length of the path is minimal in the graph ${\mathbf{G}}$.

\begin{defn}We say an exchange graph ${\mathbf{G}}$ has the {\em non-leaving-face property} if any geodesic connecting two vertices in ${\mathbf{G}}$
lies in the minimal face containing both vertices.
\end{defn}

For any two clusters $v$ and $w$ in an exchange graph ${\mathbf{G}}$, the minimal face containing $v$ and $w$ is given as $F_{v \cap w}$.
So the non-leaving-face property says that a minimal length sequence of mutations transforming the cluster $v$ into the cluster $w$ does never mutate the elements that are already common to both clusters.
While that sounds like a very natural statement, it seems surprisingly difficult to establish in general.

The non-leaving-face property of an associahedron (of type A) was first studied by D. Sleator, R. Tarjan and
W. Thurston in \cite{Sleator-Tarjan-Thurston1988} in order to find the diameter of the associahedron.
In fact, this associahedron is the exchange graph of cluster algebra of type $A.$ The non-leaving-face property of the exchange graph for cluster algebra of
type $B,C,D$ was shown by Ceballos-Pilaud \cite{Ceballos-Pilaud2016}.

\subsection{Main result}
We now describe the main object of study of this paper, the exchange graph of an unpunctured marked surface. This exchange graph, and its corresponding cluster algebra, has been introduced by Fomin, Shapiro and Thurston in \cite{Fomin-Shapiro-Thurstion2008}.
We consider a compact connected oriented 2-dimensional bordered  Riemann surface $S$ and a finite set of marked points $M$ lying on the boundary $\partial S$ of $S$ with at least one marked point on each boundary component.
The condition $M \subset \partial S$ means that we do not allow the marked surface $(S,M)$ to have punctures (note that \cite{Fomin-Shapiro-Thurstion2008}  and some of the papers we are using
are valid in the more general context of punctured surfaces).

By a curve in $(S,M)$ we mean a continuous function $\gz : [0,1] \ra S$ with
$\gz(0),\gz(1)\in M$, and a simple curve is one where $\gamma$ is injective, except possibly at the endpoints.
We always consider curves up to homotopy, and for any collection of curves we implicitly assume that their mutual intersections are minimal possible in their respective homotopy classes.
We recall some  definitions from \cite{Fomin-Shapiro-Thurstion2008}:
\begin{defn}An \emph{arc $\da$ in $(S,M)$} is a simple non-contractible curve in $(S,M)$.
The boundary of $S$ is a disjoint union of  circles, which are subdivided by the points in $M$ into boundary segments. We call an arc $\delta$ a \emph{boundary arc} if it is homotopic
to such a boundary segment.
Otherwise, $\delta$ is said to be an \emph{internal arc}. A \emph{triangulation} of $(S,M)$ is a maximal
collection $\ggz$ of arcs that do not intersect except at their endpoints.
\end{defn}

Recall that if $\tau_i$ is an internal arc in a triangluation  $\ggz$, then there exists exactly one internal arc $\tau_i'\neq\tau_i$ in $(S,M)$ such that $f_{\tau_i}(\ggz):=(\ggz\backslash\{\tau_i\})\cup\{\tau_i'\}$ is also a triangulation of $(S,M)$.
In fact, the internal arc ${\tau_i}$ is a diagonal in the quadrilateral formed by the two triangles of $\ggz$ containing ${\tau_i}$, and ${\tau'_i}$ is the other diagonal in that quadrilateral.
We denote $\tau_i'$ by $f_\ggz(\tau_i)$ and say that $f_{\tau_i}(\ggz)$ is obtained from $\ggz$ by applying a flip along $\tau_i$. Recall that the number of internal arcs in a triangulation is constant:

\begin{prop}[\cite{Fomin-Shapiro-Thurstion2008}]In each triangulation of $(S,M)$, the number
of internal arcs is
$$n = 6g + 3b + c -6$$
where $g$ is the genus of $S$, $b$ is the number of boundary components, and
$c = |M|$ is the number of marked points.
\end{prop}

The exchange graph ${\mathbf{G}}_{(S,M)}$ of the marked surface $(S,M)$ is defined as the $n-$regular graph whose clusters are the triangulations of $(S,M)$
and where two clusters are joined by an edge precisely when two triangulations are related by a flip. Then we can state our main result as follows:
\begin{thm}\label{main-theorem}Let $(S,M)$ be a marked surface without punctures. Then the exchange graph ${\mathbf{G}}_{(S,M)}$ satisfies the non-leaving-face property.
  \end{thm}

Various aspects of the exchange graph ${\mathbf{G}}_{(S,M)}$ have been studied in
\cite{Fomin-Shapiro-Thurstion2008,Brustle-Zhang2011, Brustle-Yu2015}.
The graph is finite precisely when $S$ is a disc, all other cases yield infinite graphs.

\subsection{Key lemma}
Before we prove the main result, we give in the following a key lemma to the proof of the main result. Similar as \cite{Sleator-Tarjan-Thurston1988,Ceballos-Pilaud2016,Williams2015}, we employ the notion of projection as follows:

\begin{defn}\label{proj-def} Let ${\mathbf{G}}$ be an exchange graph and $f \subset {\mathbf{G}}$  one of its faces. We say a map
 $$p_f: {\mathbf{G}} \longrightarrow f$$ is a projection if the following properties hold
 \begin{itemize}
 \item[$(p1):$]$p_f(v_i)$ is a vertex in $f$ for any vertex $v_i\in {\mathbf{G}}.$
 \item[$(p2):$]$p_f(v_i)=v_i$ if $v_i$ lies in $f.$
 \item[$(p3):$] $p_f$ sends edges in ${\mathbf{G}}$ to edges or vertices in $f$, that is,
          if $v_i-v_j$ is an edge in ${\mathbf{G}}$, then either $p_f(v_i)-p_f(v_j)$ is  an edge in $f$, or $p_f(v_i)=p_f(v_j)$ is a vertex in $f.$
 \item[$(p4):$]if $v_i-v_j$ is an edge in ${\mathbf{G}}$ such that $v_i$ belongs to $f$, then $p_f(v_j)=v_i$.
 \end{itemize}
\end{defn}

The following lemma is shown in \cite{Sleator-Tarjan-Thurston1988,Williams2015,Ceballos-Pilaud2016} for the (finite) exchange graphs studied there, but the proof applies easily to our general situation:
\begin{lem}\label{in-your-face}An exchange graph ${\mathbf{G}}$ has the
  non-leaving-face property if there exists a projection $p_f$ for each face $f$ of ${\mathbf{G}}$.
\end{lem}

We describe now in more detail the finite type situation when $(S,M)$ is a disc with $c$ marked points on the boundary. We identify $(S,M)$ with a regular polygon $P_{c}$ with $c$ vertices.
In this case, the vertices of the exchange graph ${\mathbf{G}}_{(S,M)}$ correspond to the triangulations of $P_{c}$, which are given by maximal collections of diagonals, representing the arcs in $(S,M)$. The edges of the exchange graph correspond to flips in which one diagonal is removed from a triangulation and replaced by the unique other diagonal of the thus obtained quadrilateral.
The resulting exchange graph ${\mathbf{G}}_{(S,M)}$ is the graph of the associahedron of type $A_{c-3}$. See the following the associahedron of type $A_3$ for example.

\begin{center}
\includegraphics[height=3in]{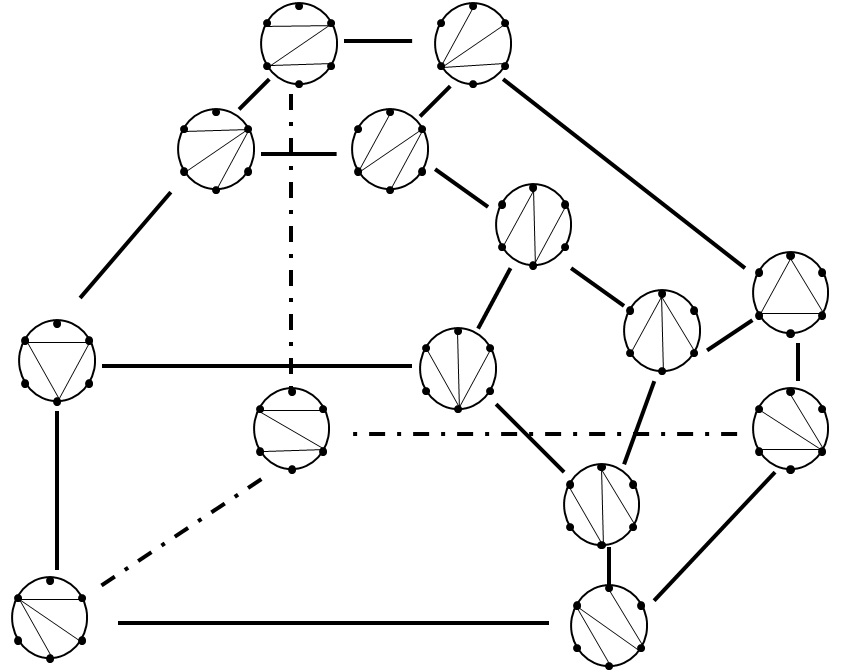}

The associahedron of type $A_3$
\end{center}

\begin{lem}[\cite{Sleator-Tarjan-Thurston1988}]\label{stt}The associahedron of type $A$ has the non-leaving-face property.
\end{lem}
Sleator, Tarjan and Thurston define in the proof of the Lemma above a projection
$p_\gz$ to a face $f$ defined by a diagonal $\gz$, that is, the face $f$ of the exchange graph that is given by all triangulations $\ggz$ containing one fixed diagonal  $\gz$.
The projection map is given in \cite{Sleator-Tarjan-Thurston1988} by some combinatorial procedure, but roughly speaking it admits the following geometric interpretation: $p_\gz(\ggz)$ is defined as the triangulation obtained by dragging all diagonals intersected by $\gz$ onto one fixed endpoint $\gz(0)$ of $\gz.$

We consider below the example of the associahedron $A_{8}$ with one fixed diagonal $\gz$ defining the face $f$, and an arbitrary triangulation $\ggz$ of $P_{11}$.
The projection $p_\gz(\ggz)$ is shown on the right side, all diagonals intersecting $\gz$ are dragged along $\gz$ onto $\gz(0)$.

\begin{center}
\includegraphics[height=2in]{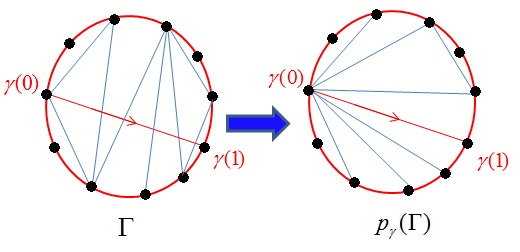}

An example of the projection
\end{center}
Inspired by the Sleator-Tarjan-Thurston's projection in the case of a disc, we define in this paper a projection for all marked surfaces $(S,M)$.

\section{Proof of the main result}
Let $(S,M)$ be a marked surface, and fix an arc $\gz$ of $(S,M)$. The face $\mathcal{F}(\gz)$ defined by $\gz$ is the full subgraph of $\mathbf{G}_{(S,M)}$ given by all triangulations that contains the arc $\gz$.
The aim of this section is to define a projection $p_\gz$ from $\mathbf{G}_{(S,M)}$ onto the face $\mathcal{F}(\gz)$.
Establishing the properties of a projection as defined in \ref{proj-def} allows us to prove the main result, by using Lemma \ref{in-your-face}.

\subsection{Projection}

We choose an orientation on the arc $\gz$, by transversing the arc from $\gz(0)$ to $\gz(1)$  in $(S,M)$.
For any arc $\gz'$ in $(S,M)$,
we denote by $\Int (\gz',\gz)$ the minimal intersection number of two representatives of the homotopy classes of $\gz'$ and $\gz$.
Moreover, for each triangulation $\ggz$ of $(S,M)$
we denote by $\t_1^\gz(\ggz)$ the first arc in $\ggz$ that intersects $\gz$ in the fixed orientation transversing $\gz$ from $\gz(0)$ to $\gz(1)$. Set
$$|\ggz|_\gz=\Int(\ggz,\gz)-\Int(\t_1^\gz(\ggz),\gz),$$
where $\Int(\ggz,\gz)=\sum_{\tau\in\ggz}\Int(\tau,\gz).$
It is known from topology that any two triangulations of the marked surface $(S,M)$ are flip-equivalent, that is one can be transformed into the other by a sequence of flips. In order to define the projection, we need an explicit proof of that fact, which will be given in the following lemma:

\begin{lem}\label{n0} Let $\ggz$ be a triangulation of $(S,M)$ that does not contain the arc  $\gz$ then there exists a sequence of flips
$$\ggz=\ggz_0\overset{f_1}{-}\ggz_1\overset{f_2}{-}\ggz_2\cdots\ggz_{m-1}\overset{f_m}{-}\ggz_m$$
such that $\gz\in\ggz_m.$
\end{lem}
\begin{pf} Note that we have $|\ggz|_\gz>0$ since $\ggz$ does not contain the arc  $\gz$.
We are first going to show  that
by an appropriate sequence of flips we obtain a triangulation $\ggz'$ with $|\ggz'|_\gz=0.$

Let us enumerate the arcs of  $\ggz_0=\ggz=\{\t_1,\t_2,\ldots,\t_n\}$ such that $\t_{1}=\t^\gz_1({\ggz_0})$ is the first arc  which intersects $\gz$ along the fixed orientation starting from $\gz(0)$, and such that $\t_{2}$ is the next arc which intersects $\gz$. The arcs $\t_{1},\t_{2}$ belong to triangles of $\ggz$ which are bordered by other arcs. We label those arcs  $\t_{3},\t_{4},\t_{5}$  in  the following figure, keeping in mind that  they might be not all distinct: depending on the surface $S$, we may have $\t_{3}=\t_{4}$ or $\t_{2}=\t_{5}$ etc.

\begin{center}
\includegraphics[height=1.8in]{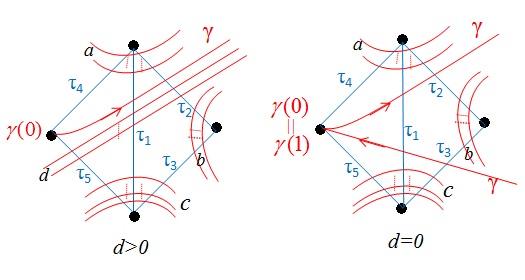}

$\gz$ intersects $\ggz_0$ with two cases at endpoints
\end{center}

In the figure, we denote by $a$ the number of times that $\gz$ intersects $\t_4,\t_1,\t_2$ successively (or $\t_2,\t_1,\t_4$ successively)
along the orientation of $\gz,$  by $d$ the number of times that $\gz$ intersects $\t_5,\t_1,\t_2$ successively (or $\t_2,\t_1,\t_5$ successively), and  similarly we define $b,c.$
Note that one might have $\gz(0)=\gz(1)$, in which case $d=0$, see the right picture.
We only consider in the following the case $\gz(0)\neq \gz(1)$, the proof for the other case is similar.

By definition, we get
$$|\ggz_0|_\gz=2a+2b+2c+2d+1+K$$
where $K=\Int(\gz,\ggz)-\sum^5_{i=1}\Int(\t_i,\ggz).$

Applying a flip on $\t_1$ we obtain a new arc $\t'_1$ and a new triangulation $\ggz_1,$ see the following picture:

\begin{center}
\includegraphics[height=2in]{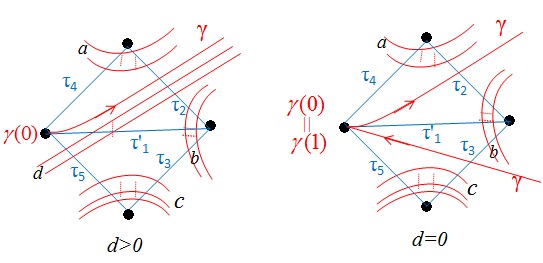}

$\gz$ intersects $\ggz_1$ with two cases at endpoints after
\end{center}

Thus $\t^\gz_1({\ggz_1})=\t_2,$ and
$$|\ggz_1|_\gz=a+2b+2c+2d+K.$$
Therefore $|\ggz_1|_\gz<|\ggz_0|_\gz.$ Since $\Int(\ggz,\gz)$ is finite, we get  by a sequence of flips, always using the next arc that crosses $\gz$, from $\ggz_0$ to a triangulation $\ggz'=\ggz_{m-1}$ with $|\ggz'|_\gz=0$, that is, $\Int(\ggz',\gz)=\Int(\t^\gz_1({\ggz'}),\gz).$
$$\ggz=\ggz_0\overset{f_{\t_1^\gz({\ggz_0})}}{-}\ggz_1\overset{f_{\t_1^\gz({\ggz_1})}}{-}\ggz_2\overset{f_{\t_1^\gz({\ggz_2})}}{-}\cdots {-} \;\ggz'=\ggz_{m-1}\overset{f_{\t_1^\gz({\ggz_{m-1}})}}{-}\ggz_m$$
Finally, one obtains a triangulation $\ggz_m$ containing $\gz$  by applying a flip from $\ggz'$ at the arc $\t^\gz_1({\ggz'})$, which completes the proof.
\end{pf}

\begin{defn}\label{proj-map} Let $\ggz$ be a triangulation of $(S,M)$, and $\gz$ an arc in $(S,M)$.
As seen in the proof of the previous lemma, applying sequences of flips at arcs intersecting $\gz$ (in the order given by the orientation of $\gz$) yields a unique triangulation $\ggz_m=p_\gz(\ggz)$ which contains the arc $\gz$. We thus obtain  a map
$$p_\gz: {\mathbf{G}}_{(S,M)}\longrightarrow \mathcal{F}(\gz)$$
where $\mathcal{F}(\gz)$ is the face associated to $\gz$.
  \end{defn}

In fact, the projection given (for the case when $S$ is a disc) in  \cite{Sleator-Tarjan-Thurston1988}  by dragging intersecting arcs along $\gz$
coincides with our definition.
We illustrate this by the following example.

\begin{example}We consider again  the example $A_8$.
Since each diagonal has two endpoints, one may have two  projections. However, if we fix an orientation
$\overset{\ra}{\gz}$ of $\gz,$ then Sleator-Tarjan-Thurston's projection can be realized by
the projection $p_{\overset{\ra}{\gz}}(\ggz)$ obtained from $\ggz$ by applying a sequence of ordered edge flips (induced by the intersections  along the orientation of $\overset{\ra}{\gz}$). We illustrate that in following figure:
\begin{center}
\includegraphics[height=2.8in]{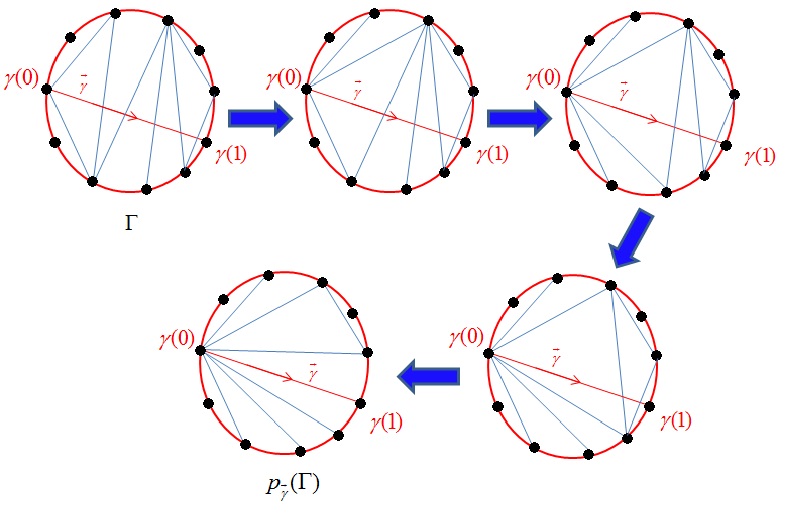}

An example of projection $p_{\overset{\ra}{\gz}}(\ggz)$
\end{center}

The thus obtained projection corresponds to dragging intersecting arcs toward the starting point of $\gz$.
However, if one chooses the opposite orientation $\overset{\leftarrow}{\gz}$ of $\gz,$ and performs the flips along the new orientation of $\overset{\leftarrow}{\gz},$ the corresponding projection $p_{\overset{\leftarrow}{\gz}}(\ggz)$ is given as follows:

\begin{center}
\includegraphics[height=2.8in]{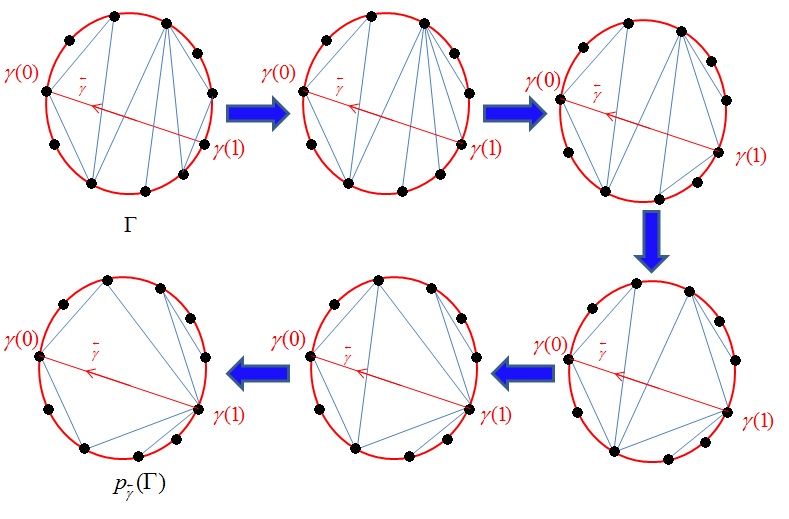}

An example of projection $p_{\overset{\leftarrow}{\gz}}(\ggz)$
\end{center}
The projection corresponds in this case to dragging intersecting arcs toward the other endpoint of $\gz$.
So, once the appropriate orientation is
fixed, the projection we defined here is the same as the one given in \cite{Sleator-Tarjan-Thurston1988}.

\end{example}

\begin{rem}Note that, once the orientation of an arc $\gz$ is chosen,  the projection $p_\gz(\ggz)$ is unique for each triangulation $\ggz.$
Moreover, given internal arcs $\gz_1,\gz_2,\ldots,\gz_m$,
 we choose for each of them an orientation, and
 denote by ${\mathcal{F}}_{(\gz_1,\gz_2,\ldots,\gz_m)}$ the minimal face of ${\mathbf{G}}_{(S,M)}$
containing the arcs $\gz_1,\gz_2,\ldots\gz_m$.
We then define a projection
$$p_{{(\gz_1,\gz_2,\ldots,\gz_m)}}=p_{\overset{\ra}{\gz_1}}\circ p_{\overset{\ra}{\gz_2}}\circ\cdots\circ p_{\overset{\ra}{\gz_m}}$$
as  composition of projections $p_{\overset{\ra}{\gz_i}}$ for all $1\leq i\leq n.$
This is a map from the exchange graph to its face ${\mathcal{F}}_{(\gz_1,\gz_2,\ldots,\gz_m)}.$

\end{rem}

\subsection{Proof of Theorem \ref{main-theorem}}
We prove our main result using Lemma \ref{in-your-face}, that is, we show that the projection map we constructed in Definition \ref{proj-map} satisfies the  properties $p(1), p(2), p(3)$ and $p(4)$ from Definition \ref{proj-def}.

In fact, properties $p(1)$ and $p(2)$ follow directly from the construction since the resulting triangulation $\ggz_m$ in  Lemma \ref{n0} contains the arc $\gz$, and  we have $\ggz_m=\ggz_0$ if $\gz$ is contained in $\ggz_0$ already.

The property $p(4)$ means the following in our context: Suppose $\ggz$ and $\ggz'$ are related by a flip at $\t_1\in\ggz$, and assume further that the arc $\gz$ belongs to $\ggz'$, but not to $\ggz$ (if both $\ggz$ and $\ggz'$ are already in the face ${\mathcal{F}}_{\gz}$ there is nothing to show). But that means that the flip at the first arc $\t_1$ creates $\gz$, so by our construction of the projection map we have $p_\gamma(\ggz) = \ggz'$, which was to show.

Therefore it suffices to prove property $(p3).$
We consider the projection $p_\gamma$ with respect to a fixed arc $\gz$ as in  Lemma \ref{n0}, and assume that two triangulations $\ggz$ and $\ggz'$ are related by a flip at some arc $\t_1\in\ggz$. That is, $\gz$
intersects $\ggz$ and $\ggz'$ in the same way except at $\t_1\in\ggz$ and $\t_1'\in\ggz'\setminus\ggz$, see the following:
\begin{center}
\includegraphics[height=1.8in]{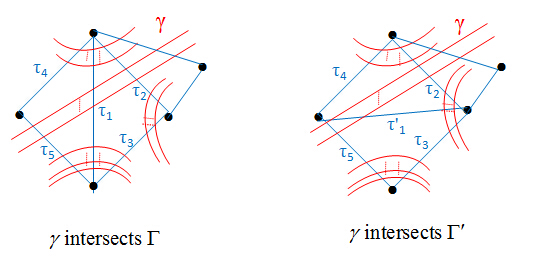}

$\gz$ intersects triangulations related by a flip
\end{center}
By construction of the projection $p_\gamma$, one proceeds performing the  same sequence of flips along the orientation of $\gz$ to both triangulations  $\ggz$ and $\ggz'$ as long the arc $\t_1\in\ggz$ and $\t_1'\in\ggz'$ is not yet encountered:
 $$\ggz\overset{flips}{--}\ggz_i,\\\ \ggz'\overset{flips}{--}\ggz'_i.$$
Therefore the intermediate triangulations $\ggz_i$ and $\ggz'_i$ are still related by a flip at $\t_1$, and we only need to consider the  instance when the construction of  $p_\gamma$ flips at the arc $\t_1\in\ggz$ or $\t_1'\in\ggz'$, respectively.
It is thus sufficient to consider one of the following situations (or their duals):
\bigskip

{\bf Case I:}
\begin{center}
\includegraphics[height=1.8in]{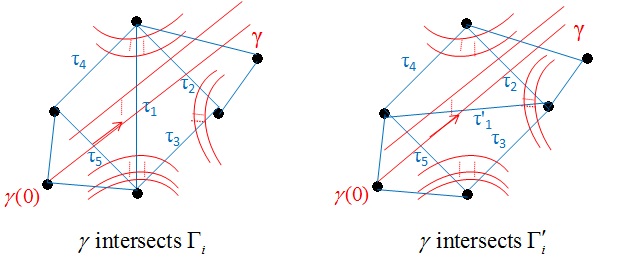}

After some flips, $\gz$ intersects $\ggz'_i$ at $\t_5,\t_1...$ and intersects $\ggz'_i$ at $\t_5,\t'_1...$
\end{center}
By definition of the projection, one should apply successively flips at $\t_5,\t_1,\t_2,\ldots$ for $p_\gz(\ggz):$
$$\ggz\overset{flips}{--}\ggz_i\overset{f_{\t_5}}{-}\ggz_{i+1}\overset{f_{\t_1}}{-}\ggz_{i+2}\overset{f_{\t_2}}{-}\ggz_{i+3}\ldots$$
and apply successively flips at $\t_5,\t'_1,\t_2,\ldots$ for $p_\gz(\ggz'):$
$$\ggz'\overset{flips}{--}\ggz'_i\overset{f_{\t_5}}{-}\ggz'_{i+1}\overset{f_{\t'_1}}{-}\ggz'_{i+2}\overset{f_{\t_2}}{-}\ggz'_{i+3}\ldots.$$
But then it is easy to verify that $\ggz_{i+3}=\ggz'_{i+3}$, which implies that $p_\gz(\ggz)=p_\gz(\ggz')$.
\bigskip

{\bf Case II:}
\begin{center}
\includegraphics[height=1.8in]{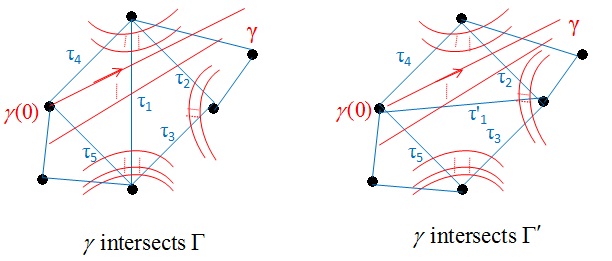}

$\t^\ggz_1(\gz)=\t_1$ but $\t^{\ggz'}_1(\gz)\neq\t'_1$
\end{center}
In this situation, by definition of the projection one performs the following flips for $p_\gz(\ggz)$:
$$\ggz\overset{f_{t_1}}{-}\ggz_1\overset{f_{\t_2}}{-}\ggz_{2}\ldots$$
and likewise for $p_\gz(\ggz')$:
$$\ggz'\overset{f_{\t_2}}{-}\ggz'_1\ldots$$
Hence $\ggz'_1=\ggz_2$ which yields $p_\gz(\ggz)=p_\gz(\ggz')$ by definition of the projection.
\bigskip

{\bf Case III:}
\begin{center}
\includegraphics[height=1.8in]{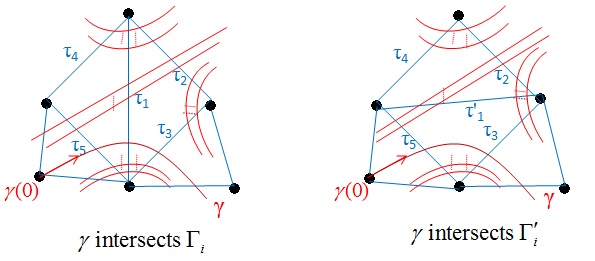}
\end{center}

Similarly as Case I, by definition of the projection, one should apply successively flips at $\t_5,\t_1,\t_3,\ldots$ for $p_\gz(\ggz):$
$$\ggz\overset{flips}{--}\ggz_i\overset{f_{\t_5}}{-}\ggz_{i+1}\overset{f_{\t_1}}{-}\ggz_{i+2}\overset{f_{\t_3}}{-}\ggz_{i+3}\ldots$$
and apply successively flips at $\t_5,\t_3,\ldots$ for $p_\gz(\ggz'):$
$$\ggz'\overset{flips}{--}\ggz'_i\overset{f_{\t_5}}{-}\ggz'_{i+1}\overset{f_{\t_3}}{-}\ggz'_{i+2}\ldots.$$ After these flips, we get the following pictures for $\ggz_{i+3}$ and $\ggz'_{i+2}$ respectively: see
\begin{center}
\includegraphics[height=1.8in]{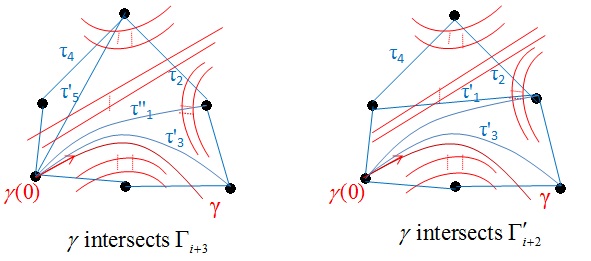}
\end{center}
Then $\ggz_{i+3}$ and $\ggz'_{i+2}$ are related by a flip at $\t'_1\in\ggz'_{i+2}$ or $\t'_5\in\ggz_{i+3}$ and we can proceed by induction on the number of flips needed in the construction of $p_\gz(\ggz)$ to show the statement of property $p(3)$. This completes the proof of Theorem \ref{main-theorem}.



\end{document}